\newcount\ite\ite=1\def\0{\global\ite=1\1}
\def\1{\item{\rm(\romannumeral\the\ite)}\advance\ite1\quad}
\def\phi{\varphi}

\documentclass[12pt,twoside]{article}
\usepackage{amssymb}
\usepackage{setspace}

\font\teneufm=eufm10 scaled \magstep1
\font\seveneufm=eufm7 scaled \magstep1
\font\fiveeufm=eufm5  scaled \magstep1
\newfam\eufmfam

\textfont\eufmfam=\teneufm
\scriptfont\eufmfam=\seveneufm
\scriptscriptfont\eufmfam=\fiveeufm

\newfam\msbfam
\font\tenmsb=msbm10 scaled \magstep1  \textfont\msbfam=\tenmsb
\font\sevenmsb=msbm7 scaled \magstep1 \scriptfont\msbfam=\sevenmsb
\font\fivemsb=msbm5 scaled \magstep1  \scriptscriptfont\msbfam=\fivemsb

\def\dd#1{\raise1.5pt\hbox{$\,\partial\!$}/\raise-2.5pt\hbox{$\!\partial#1\,$}}

\def\tilde{\widetilde}
\def\hat{\widehat}
\def\5#1{{\mathcal #1}}

\def\CC{{\mathbb C}}

\def\ra{\rightarrow}

\def\GL{\mathop{\rm GL}\nolimits}

\def\Ann{\mathop{\rm Ann}\nolimits}

 \def\HollowBoxx #1#2#3{{\dimen0=#1 \advance\dimen0 by -#2
       \dimen1=#1 \advance\dimen1 by #3
        \vrule height 0pt depth #3 width #2
       \hskip -#3
       \vrule height #1 depth #3 width #3}}
 \def\LeftContraction{\mathord{\kern1.45pt \HollowBoxx{6pt}{3.5pt}{.4pt}}\,}

 \def\HollowBox #1#2#3{{\dimen0=#1 \advance\dimen0 by -#3
       \dimen1=#1 \advance\dimen1 by #3
        \vrule height #1 depth #3 width #3
        \vrule height 0pt depth #3 width #2
        \hskip -#3}}
 \def\RightContraction{\mathord{\, \HollowBox{6pt}{3.1pt}{.4pt}} \kern1.6pt}

\def\qed{{\hfill $\Box$}}
\newtheorem{theorem}{THEOREM}[section]

\newtheorem{remark}[theorem]{Remark}
\newtheorem{proposition}[theorem]{Proposition}
\newtheorem{conjecture}[theorem]{Conjecture}

\begin{document}

\begin{center}
{\Large \bf Application of Classical Invariant Theory\\
\vspace{0.1cm}
to Biholomorphic Classification\\
\vspace{0.1cm}
of Plane Curve Singularities,\\
\vspace{0.3cm}
and Associated Binary Forms}
\vspace{0.2cm}\\
\normalsize A. V. Isaev\footnote{
\noindent Department of Mathematics, The Australian National University, Canberra, ACT 0200, Australia. 

E-mail: alexander.isaev@anu.edu.au}
\end{center}

\begin{quotation} 
{\small \sl \noindent We use classical invariant theory to solve the biholomorphic equivalence problem for two families of plane curve singularities previously considered in the literature. Our calculations motivate an intriguing conjecture that proposes a way of extracting a complete set of invariants of homogeneous plane curve singularities from their moduli algebras.}
\end{quotation}

\thispagestyle{empty}

\pagestyle{myheadings}
\markboth{A. V. Isaev}{Plane Curve Singularities and Associated Binary Forms}

\setcounter{section}{0}

\section{Introduction}\label{intro}
\setcounter{equation}{0}

In singularity theory one often encounters parametric families of singularities (see, e.g. \cite{A}), and for a given family ${\mathcal F}$ it is desirable to characterize all subsets of parameter values for which the corresponding singularities in ${\mathcal F}$ are pairwise biholomorphically equivalent. This problem may be quite non-trivial even if ${\mathcal F}$ is defined by equations of simple form. In the present paper we consider families of homogeneous plane curve singularities and discuss the equivalence problem for such families from the point of view of classical invariant theory.

Recall that a {\it plane curve}\, is a complex curve in $\CC^2$, i.e. a subset of the form $V=\{p\in\ U: f(p)=0\}$, where $f\not\equiv 0$ is a function holomorphic in a domain $U\subset\CC^2$. Everywhere in this paper plane curves will be considered with the corresponding reduced complex structures. We assume that $V$ contains the origin and will only be interested in the germ ${\mathcal V}$ of $V$ at 0. Shrinking $U$ if necessary one can assume that the defining function $f$ of $V$ is {\it minimal}\, in the sense that the set
$
C_f:=\left\{p\in V: df(p)=0\right\}
$  
is nowhere dense in $V$, in which case $C_f$ coincides with the singular set of $V$. If $f$ and $g$ are two minimal defining functions of $V$, then $f=u\cdot g$, where $u$ is a holomorphic nowhere vanishing function in $U$. The germ ${\mathcal V}$ is called {\it singular}\, if $0\in C_f$.

Next, two plane curve germs ${\mathcal V}_1$ and ${\mathcal V}_2$ are said to be {\it biholomorphically equivalent}\, if there exists a biholomorphic mapping $F$ between neighborhoods $U_1$ and $U_2$ of the origin in $\CC^2$, such that $F(0)=0$ and $F(V_1)= V_2$, where $V_1\subset U_1$ and $V_2\subset U_2$ are plane curves representing ${\mathcal V}_1$ and ${\mathcal V}_2$, respectively. We are interested in the biholomorphic equivalence problem for singular plane curve germs.

Let ${\mathcal O}_2$ be the local algebra of holomorphic function germs at the origin in $\CC^2$. For a plane curve germ ${\mathcal V}$ as above, denote by ${\mathcal A}({\mathcal V})$ the quotient of ${\mathcal O}_2$ by the ideal generated by the germs of $f$ and all its first-order partial derivatives. The algebra ${\mathcal A}({\mathcal V})$, called the {\it moduli algebra}\, or {\it Tjurina algebra}\, of ${\mathcal V}$, is independent of the choice of $f$ as well as the coordinate system near the origin, and the moduli algebras of biholomorphically equivalent curve germs are isomorphic. Clearly, ${\mathcal A}({\mathcal V})$ is non-trivial if and only if the germ ${\mathcal V}$ is singular. Furthermore, since the singularity of ${\mathcal V}$ is necessarily isolated, we have $\dim_{\CC}{\mathcal A}({\mathcal V})<\infty$ (see, e.g. \cite{GLS}). By a theorem due to Mather and Yau (see \cite{MY}), two germs ${\mathcal V}_1$, ${\mathcal V}_2$ are biholomorphically equivalent if their moduli algebras ${\mathcal A}({\mathcal V})$, ${\mathcal A}({\mathcal V}_2)$ are isomorphic. Thus, the moduli algebra ${\mathcal A}({\mathcal V})$ determines ${\mathcal V}$ up to biholomorphism. The proof of the Mather-Yau theorem does not provide an explicit procedure for recovering the germ ${\mathcal V}$ from the algebra ${\mathcal A}({\mathcal V})$, and finding such a procedure is an interesting open problem. Calculations performed in this paper for particular families of singularities motivate a conjecture that proposes a procedure of this kind in the case of homogeneous plane curve germs.

A germ ${\mathcal V}$ is said to be {\it homogeneous}\, if it is represented by a plane curve $V$ that in some holomorphic coordinates $z$, $w$ near 0 is defined by a homogeneous polynomial (binary form) $Q$. Any non-zero binary form can be written as a product of linear factors, and the minimality of $Q$ means that each of its factors has multiplicity one, i.e. $Q$ is {\it square-free}. Clearly, if $\deg Q\ge 2$, the germ ${\mathcal V}$ is singular. If ${\mathcal V}_1$, ${\mathcal V}_2$ are two plane curve germs defined by square-free binary forms $Q_1$, $Q_2$, respectively, then the germs are biholomorphically equivalent if and only if the binary forms are {\it linearly equivalent}, that is, there exists $C\in\GL(2,\CC)$ such that $Q_1(C{\bf z})\equiv Q_2({\bf z})$, with ${\bf z}:=(z,w)$. The linear equivalence problem for non-zero square-free binary forms can be solved by utilizing absolute classical invariants (see Proposition \ref{equivarbitr} for details). Thus, the biholomorphic equivalence problem for homogeneous plane curve germs is solvable, in principle, by methods of classical invariant theory.

To demonstrate the power of the classical invariant theory approach, we start with the family of plane curves given by
\begin{equation}
f_{t}:=z^n+t z^{n-1}w+w^n=0,\quad t\in\CC,\quad n\ge 4. \label{family}
\end{equation}
Here the form $f_{t}$ is assumed to be square-free, which in terms of the parameter $t$ means $t^n\ne-n^n/(1-n)^{n-1}$ (see (\ref{denominator})). Family (\ref{family}) first appeared in paper \cite{KK}, where the plane curves $\{z^n+a(w)w^{\alpha}z+b(w)w^{\beta}=0\}$ were considered. Here $\alpha$, $\beta$ are positive integers, $n\ge 3$, and $a(w)$, $b(w)$ are holomorphic nowhere vanishing functions defined near the origin. The germs of the above curves at the origin were classified in \cite{KK} up to biholomorphic equivalence in many situations (see also \cite{SL}). However, the homogeneous case $\alpha=n-1$, $\beta=n$, $a\equiv 1$, $b\equiv\hbox{const}$, where the form $z^n+w^{n-1}z+bw^n$ is square-free, proved to be of substantial difficulty to the authors for $n\ge 4$ (observe that for $n=3$ the germs defined by these forms are pairwise biholomorphically equivalent for trivial reasons). The biholomorphic equivalence problem in the homogeneous case was eventually solved in later paper \cite{Ka} as stated in the theorem below (note that for $n\ge 4$ and $b\ne 0$ the binary form $z^n+w^{n-1}z+bw^n$ is linearly equivalent to $f_t$ for some $t$). 

\begin{theorem}\label{main1}\sl Let ${\mathcal V}_t$ be the germ of curve {\rm (\ref{family})} at the origin.  Then two germs ${\mathcal V}_{t_1}$ and ${\mathcal V}_{t_2}$ are biholomorphically equivalent if and only if $t_1^n=t_2^n$.
\end{theorem}

Theorem \ref{main1} was obtained in \cite{Ka} by the direct substitution method, which required enormous calculations. In Section \ref{proofs1} below we recover Theorem \ref{main1} by a short elementary argument based on classical invariant theory. Our proof is preceded by a brief survey of some basic facts concerning classical invariants in Section \ref{prelim}.

If one introduces more terms in equation (\ref{family}), the biholomorphic equivalence problem becomes much harder even for low values of $n$. Indeed, let ${\mathcal V}_{s,t}$ be the germ at the origin of the curve
\begin{equation}
f_{s,t}:=z^5+sz^4w+tz^3w^2+w^5=0,\quad s,t\in\CC,\label{familyfst}
\end{equation}
where $f_{s,t}$ is assumed to be square-free. This family was considered in \cite{Ea} and \cite{Ka}. As observed in \cite{Ka}, the germs ${\mathcal V}_{5,10}$ and ${\mathcal V}_{15\cdot 5^{-4/5},10\cdot 5^{-3/5}}$ are biholomorphically equivalent, which indicates that for two equivalent germs ${\mathcal V}_{s_1,t_1}$ and ${\mathcal V}_{s_2,t_2}$ the pairs of parameter values $s_1$, $t_1$ and $s_2$, $t_2$ may be related in a complicated way. We discuss the family ${\mathcal V}_{s,t}$ in Section \ref{n5} and show that in this case classical invariant theory provides only an implicit solution to the biholomorphic equivalence problem, with the pairs $s_1$, $t_1$ and $s_2$, $t_2$ related by means of a system of three rational equations (see Theorem \ref{specialequivsol}). Our consideration of the family ${\mathcal V}_{s,t}$ in Section \ref{n5} is preceded by a discussion of classical results on the linear equivalence problem for general square-free binary quintics in Section \ref{quinticsss}, where the main fact is Theorem \ref{solinv}.

Section \ref{assocforms} is perhaps the most interesting part of the paper. In this section, to every non-zero square-free binary form $Q$ of degree $n$ we associate a binary form $\hat Q$ of degree $2(n-2)$ defined up to scale. The form $\hat Q$ arises from the moduli algebra of the germ at the origin of the plane curve defined by $Q$. Such associated forms, introduced in a slightly different way, were first considered by M. Eastwood in \cite{Ea} for the purpose of finding an explicit procedure for extracting invariants of isolated hypersurface singularities from their moduli algebras. In particular, in \cite{Ea} the sextics $\hat f_{s,t}$ associated to the quintics $f_{s,t}$ were computed, and by applying classical invariant theory to $\hat f_{s,t}$ Eastwood was able to recover some invariants of $f_{s,t}$. In Section \ref{assocforms} we show that, interestingly, one in fact can obtain all three rational equations of Theorem \ref{specialequivsol} by supplementing the computations performed in \cite{Ea}. Further, in this section we also give an alternative proof of Theorem \ref{main1} by applying classical invariant theory to the forms $\hat f_t$ associated to $f_t$. Inspired by the examples of the families ${\mathcal V}_t$, ${\mathcal V}_{s,t}$, we propose a general conjecture relating the invariant theory of non-zero square-free binary forms to that of the corresponding associated forms (see Conjecture \ref{conj}). If correct, the conjecture would yield an explicit algorithm for extracting a complete set of invariants of homogeneous plane curve germs from their moduli algebras, which would complement the Mather-Yau theorem in this case. In our forthcoming joint paper with M. Eastwood, Conjecture \ref{conj} will be shown to hold for binary quintics and binary sextics. Note that an analogous conjecture can be stated in any number of variables.

{\bf Acknowledgements.} We would like to thank M. Eastwood for many discussions of results in \cite{Ea} as well as for performing some Maple calculations (see Section \ref{assocforms}), A. Gorinov for a useful conversation (see the proof of Proposition \ref{equivarbitr}), and N. Kruzhilin for bringing papers \cite{KK}, \cite{SL} to our attention. This work is supported by the Australian Research Council.

\section{Preliminaries}\label{prelim}
\setcounter{equation}{0}

In this section we collect some basic facts from classical invariant theory. Further details can be found, for example, in monograph \cite{O}. Everywhere below the ground field is assumed to be $\CC$. 

We consider {\it binary forms}, i.e. homogeneous polynomials on a 2-dimen\-sion\-al vector space $V$. Let ${\mathcal Q}_V^n$ be the linear space of binary forms of a fixed degree $n\ge 2$. Define an action of $\GL(V)$ on ${\mathcal Q}_V^n$ by the formula
$$
(C,Q)\mapsto Q_C,\quad  Q_C(v):=Q(C^{-1}v),\,\,C\in\GL(V),\,\, Q\in{\mathcal Q}_V^n,\,\, v\in V.
$$
Two binary forms are said to be {\it linearly equivalent}\, if they lie in the same orbit with respect to this action.

An {\it invariant}\, (or {\it relative classical invariant}) of binary forms of degree $n$ on $V$ is a polynomial $I:{\mathcal Q}_V^n\ra\CC$ such that for any $Q\in{\mathcal Q}_V^n$ and any $C\in\GL(V)$ one has $I(Q)=(\det C)^kI(Q_C)$, where $k$ is a non-negative integer called the {\it weight}\, of $I$. It follows that $I$ is in fact homogeneous of degree $2k/n$.

Next, for any two invariants $I$ and $\tilde I$ the ratio $I/\tilde I$ yields a rational function on ${\mathcal Q}_V^n$ that is defined, in particular, at the points where $\tilde I$ does not vanish. If $I$ and $\tilde I$ have equal weights, this function does not change under the action of $\GL(V)$, and we say that $I/\tilde I$ is an {\it absolute invariant}\, (or {\it absolute classical invariant}) of forms of degree $n$ on $V$. If one fixes coordinates $z$, $w$ in $V$, then any element $Q\in{\mathcal Q}^n_V$ is written as
$$
Q(z,w)=\sum_{i=0}^n\left(\begin{array}{l} n\\ i\end{array}\right) a_iz^iw^{n-i},
$$
where $a_i\in\CC$. In what follows we will introduce a number of absolute invariants of binary forms defined in terms of the coefficients $a_i$. Observe that for any absolute invariant ${\mathcal I}$ so defined its value ${\mathcal I}(Q)$ is in fact independent of the choice of coordinates in $V$. When working in coordinates, we always assume that $V=\CC^2$ and identify $\GL(V)$ with $\GL(2,\CC)$.    

Any form $Q\in{\mathcal Q}^n_{\CC^2}$ can be written as a product of linear terms
\begin{equation}
Q(z,w)=\prod_{\alpha=1}^n(w_{\alpha}z-z_{\alpha}w),\label{representation}
\end{equation}
for some $z_{\alpha},w_{\alpha}\in\CC$. The {\it discriminant}\, of $Q$ is then defined as
$$
\Delta(Q):=\frac{(-1)^{n(n-1)/2}}{n^n}\prod_{1\le\alpha<\beta\le n}(z_{\alpha}w_{\beta}-z_{\beta}w_{\alpha})^2
$$ 
(see pp. 97--101 in \cite{El}). Clearly, this definition is independent of representation (\ref{representation}). The discriminant $\Delta(Q)$ is an invariant of degree $2(n-1)$, which is non-zero if and only if $Q$ is non-zero and square-free. Furthermore, if $a_n\ne 0$, the discriminant $\Delta(Q)$ can be computed as
\begin{equation}
\displaystyle\Delta(Q)=\frac{{\bf R}(Q,\partial Q/\partial z)}{n^n a_n},\label{comp}
\end{equation}    
where for two forms $P(z,w)$ and $S(z,w)$ we denote by ${\bf R}(P,S)$ their resultant (see p. 36 in \cite{O}). 

Next, we define the $n$th {\it transvectant}\, as
\begin{equation}
(Q,Q)^{(n)}:=(n!)^2\sum_{i=0}^n(-1)^i\left(\begin{array}{c} n \\ i\end{array}\right)a_ia_{n-i}.\label{transvectant}
\end{equation}
The transvectant $(Q,Q)^{(n)}$ is an invariant of degree 2. It is identically zero if $n$ is odd, thus for any odd $n$ we consider the invariant $\left(Q^2,Q^2\right)^{(2n)}$, which has degree 4. 

We now introduce an absolute invariant as follows:
\begin{equation}
\displaystyle J(Q):=\left\{\begin{array}{ll}\displaystyle\frac{\Bigl[(Q,Q)^{(n)}\Bigr]^{n-1}}{\Delta(Q)} & \hbox{if $n$ is even,}\\
\vspace{-0.3cm}\\
\displaystyle\frac{\Bigl[\left(Q^2,Q^2\right)^{(2n)}\Bigr]^{(n-1)/2}}{\Delta(Q)} & \hbox{if $n$ is odd.}
\end{array}\right.\label{invariantj}
\end{equation}
Next, for even values of $n$ we introduce the absolute invariant
\begin{equation}
M(Q):=\frac{({\mathbf H}(Q),{\mathbf H}(Q))^{(2(n-2))}}{\Bigl[(Q,Q)^{(n)}\Bigr]^2},\label{invariantm}
\end{equation}
where ${\bf H}(Q)$ is the Hessian of $Q$. Note that ${\mathbf H}(Q)\in{\mathcal Q}^{2(n-2)}_{\CC^2}$ and the relative invariant $({\mathbf H}(Q),{\mathbf H}(Q))^{(2(n-2))}$ has degree 4.

In our proof of Theorem \ref{main1} in the next section we compute $J(f_t)$, where $f_t$ is the binary form defined in (\ref{family}). In Section \ref{assocforms} we will give an alternative proof of Theorem \ref{main1} based on computing $M({\bf f}_t)$, where ${\bf f}_t$ is a certain form of degree $2(n-2)$ arising from $f_t$. 

\section{Proof of Theorem \ref{main1}}\label{proofs1}
\setcounter{equation}{0}

The \lq\lq if\rq\rq\, implication of the theorem is trivial since the curve $\{f_{\rho t}(z,w)=0\}$, with $\rho^n=1$, is biholomorphically equivalent to the curve $\{f_{t}(z,w)=0\}$ by means of the map $z\mapsto \rho z$, $w\mapsto w$.

To obtain the \lq\lq only if\rq\rq\, implication, we find $J(f_t)$. A straightforward computation yields
\begin{equation}
\displaystyle \begin{array}{ll} (f_t,f_t)^{(n)}=2(n!)^2 &  \hbox{if $n$ is even},\\
\vspace{-0.3cm}\\
\displaystyle (f_t^2,f_t^2)^{(2n)}=2 (2n)!\Bigl((2n)!-2(n!)^2\Bigr) &  \hbox{if $n$ is odd}.
\end{array}\label{numerator}
\end{equation}
Therefore, the numerators in (\ref{invariantj}) do not depend on $t$ and are non-zero. We will now compute the discriminant $\Delta(f_t)$. Since for $f_t$ we have $a_n=1$, one can apply formula (\ref{comp}). The resultant ${\bf R}(f_t,\partial f_t/\partial z)$ can be easily found by using cofactor expansions, and we get ${\bf R}(f_t,\partial f_t/\partial z)=(1-n)^{n-1}t^n+n^n$.  Hence
\begin{equation}
\Delta(f_t)=(1-n)^{n-1}t^n/n^n+1.\label{denominator}
\end{equation} 
Formulas (\ref{numerator}) and (\ref{denominator}) imply
\begin{equation}
J(f_t)=\frac{1}{\mu t^n+\nu}\label{jphit}
\end{equation}
for some $\mu,\nu\in\CC$ with $\mu\ne 0$, which yields the desired solution to the biholomorphic equivalence problem for the germs ${\mathcal V}_t$.\qed

\section{Equivalence of square-free binary quintics}\label{quinticsss}
\setcounter{equation}{0}

In this section we discuss the problem of linear equivalence for non-zero square-free binary quintics and explain how this problem can be solved by means of specific absolute classical invariants. What follows is known to the experts but is not easy to find in the literature.

We start with a general proposition that applies to binary forms of an arbitrary degree.

\begin{proposition}\label{equivarbitr} \sl Let
\begin{equation}
X^n:=\{Q\in{\mathcal Q}_{\CC^2}^n:\Delta(Q)\ne 0\}.\label{setx}
\end{equation}
For $n\ge 3$ the orbits of the $\GL(2,\CC)$-action on $X^n$ are separated by absolute classical invariants of the kind
\begin{equation}
{\mathcal I}=\frac{I}{\Delta^{m}},\label{formmsss}
\end{equation} 
where $m$ is a non-negative integer and $I$ is a relative classical invariant.
\end{proposition}

\noindent {\bf Proof:}\footnote{This proof was suggested to us by A. Gorinov.} The set $X^n$ admits the structure of an affine algebraic variety, and with respect to this structure the action of the complex reductive group $G:=\GL(2,\CC)$ on $X^n$ is algebraic. The stabilizers of the $G$-action are finite (see \cite{OS}), and therefore the quotient of $X^n$ by this action coincides with the Hilbert quotient $Z:=X^n/\hspace{-0.1cm}/G$ (see, e.g. \cite{Kr}, \cite{MFK}). On the quotient $Z$ one can introduce the structure of an affine algebraic variety in such a way that the quotient map $\pi: X^n\ra Z$ is an algebraic morphism and $\pi^*: \CC[Z]\ra \CC[X^n]^{G}$ is an isomorphism, where $\CC[Z]$ is the algebra of regular functions on $Z$ and $\CC[X^n]^{G}$ is the algebra of $G$-invariant regular functions on $X^n$. Since the points of $Z$ are separated by elements of $\CC[Z]$, the $G$-orbits in $X^n$ are separated by elements of $\CC[X^n]^{G}$. When $X^n$ is embedded into ${\mathcal Q}_{\CC^2}^n$ as in (\ref{setx}), every element of $\CC[X^n]^{G}$ becomes the restriction to $X^n$ of an absolute invariant of the form (\ref{formmsss}). The proof is complete. \qed
\vspace{0.1cm}\\

\noindent In what follows the algebra of the restrictions to $X^n$ of absolute invariants of the form (\ref{formmsss}) is denoted by ${\mathcal I}^n$.

We will now concentrate on the case of binary quintics (here $n=5$). Define the {\it canonizant} of binary quintics as
$$
\hbox{Can}(Q):=\det\left(\begin{array}{lll}
a_5z+a_4w & a_4z+a_3w & a_3z+a_2w\\
a_4z+a_3w & a_3z+a_2w & a_2z+a_1w\\
a_3z+a_2w & a_2z+a_1w & a_1z+a_0w
\end{array}\right)
$$
for all $Q\in{\mathcal Q}_{\CC^2}^5$ (see p. 274 in \cite{El}). Clearly, $\hbox{Can}(Q)$ is a binary cubic, and we let $I_{12}(Q):=-27\cdot\Delta\left(\hbox{Can}(Q)\right)$ (see p. 307 in \cite{El}). It turns out that $I_{12}$ is a relative invariant of binary quintics of degree 12 (which explains the notation). We now introduce two absolute invariants of binary quintics
\begin{equation}
K(Q):=\frac{I_{12}(Q)^2}{\Delta(Q)^3},\quad L(Q):=\frac{\left(Q^2,Q^2\right)^{(10)}\cdot I_{12}(Q)}{\Delta(Q)^2}. \label{invariantk}
\end{equation}
We will need the following fact.

\begin{theorem}\label{solinv} \sl The algebra ${\mathcal I}^5$ is generated by the restrictions of $J$, $K$, $L$ to $X^5$. 
\end{theorem}

\noindent {\bf Proof:} Let $I_{18}$ be the classical invariant of binary quintics of degree 18 defined on p. 309 in \cite{El}. We have
\begin{equation}
16 I_{18}^2=I_4I_8^4+8I_8^3I_{12}-2I_4^2I_8^2I_{12}-72I_4I_8I_{12}^2-432I_{12}^3+I_4^3I_{12}^2 \label{i18}
\end{equation}
(see p. 313 in \cite{El}), where
$$
\begin{array}{l}
\displaystyle I_4(Q):=\frac{\left(Q^2,Q^2\right)^{(10)}}{7200000\cdot 10!},\\
\vspace{-0.1cm}\\
\displaystyle I_8(Q):=\frac{I_4(Q)^2-\Delta(Q)}{128}.
\end{array}
$$
It is well-known that the invariants $I_4$, $I_8$, $I_{12}$, $I_{18}$ generate the algebra of classical invariants of binary quintics (see \cite{Sy}). Identity (\ref{i18}) then implies that for every absolute classical invariant ${\mathcal I}$ of binary quintics of the form (\ref{formmsss}) the numerator $I$ is a polynomial in $I_4$, $I_8$, $I_{12}$, and hence the restriction ${\mathcal I}|_{X^5}$ is a polynomial in $J|_{X^5}$, $K|_{X^5}$ and $L|_{X^5}$. \qed
\vspace{0.1cm}

It is not hard to give examples showing that none of the three pairs of $J|_{X^5}$, $K|_{X^5}$, $L|_{X^5}$ generates ${\mathcal I}^5$. Indeed, for the quintics $f_t$ with $n=5$ and $256 t^5+3125\ne 0$ (see (\ref{family}), (\ref{denominator})), the invariants $K$ and $L$ vanish, but $J(f_t)$ is a non-constant function of $t$ (see (\ref{jphit})). Next, for the quintics $g_t:=z^5+5tz^4w+5zw^4/t+w^5$ with $t\ne 0$ and $t^5\ne 7\pm 4\sqrt{3}$, the invariants $J$ and $L$ vanish, but $K(g_t)$ is a non-constant function of $t$. Finally, for the quintics $h_t:=z^5/t+w^5/(1-t)+(z+w)^5$ with $t\ne 0,1, (1\pm i\sqrt{3})/2$, the value of $J$ is independent of $t$, and one can find $t_1$, $t_2$ such that $K(h_{t_1})=K(h_{t_2})$, but $L(h_{t_1})=-L(h_{t_2})\ne 0$.

\section{The family ${\mathcal V}_{s,t}$}\label{n5}
\setcounter{equation}{0}

By Proposition \ref{equivarbitr} and Theorem \ref{solinv}, to establish whether two quintics $f_{s_1,t_1}$ and $f_{s_2,t_2}$ (see (\ref{familyfst})) are linearly equivalent it is sufficient to compare the values of the absolute invariants $J$, $K$, $L$ for these forms. We will now explicitly compute $J(f_{s,t})$, $K(f_{s,t})$, $L(f_{s,t})$ for all $s,t$. 

For the numerator and denominator in formula (\ref{invariantj}) we have, respectively,
\begin{equation}
(f_{s,t}^2,f_{s,t}^2)^{(10)}=57600\cdot 10!(125-3st^2)\label{numerator5}
\end{equation}
and
\begin{equation}
\displaystyle\Delta(f_{s,t})=\frac{1}{3125}\left(256 s^5-1600s^3t-27s^2t^4+2250st^2+108t^5+3125\right).\label{denominator5}
\end{equation}
Formulas (\ref{invariantj}), (\ref{numerator5}), (\ref{denominator5}) yield
\begin{equation}
\begin{array}{l}
\displaystyle j(s,t):=J(f_{s,t})=5(1440000\cdot 10!)^2(125-3st^2)^2/\\
\vspace{-0.3cm}\\
\hspace{1.7cm}\left(256s^5-1600s^3t-27s^2t^4+2250st^2+108t^5+3125\right).
\end{array}\label{jspecial}
\end{equation}
Further, by a straightforward computation we obtain
\begin{equation}
\begin{array}{l}
I_{12}(f_{s,t})=-\displaystyle\frac{1}{10^{10}}\left(19200 s^6t^2-160000 s^4t^3-1120 s^3t^6+\right.\\
\vspace{-0.3cm}\\
\hspace{2.3cm}
\left.440000 s^2t^4+3600 st^7+27 t^{10}-400000 t^5\right).
\end{array}\label{i12}
\end{equation}
Formulas (\ref{invariantk}), (\ref{denominator5}), (\ref{i12}) imply
\begin{equation}
\begin{array}{l}
\displaystyle k(s,t):=K(f_{s,t})=\displaystyle\frac{1}{2^{20}5^5}\left(19200 s^6t^2-160000 s^4t^3-\right.\\
\vspace{-0.3cm}\\
\hspace{1.7cm}\left.1120 s^3t^6+440000s^2t^4+3600 st^7+27 t^{10}-400000 t^5\right)^2/\\
\vspace{-0.3cm}\\
\hspace{1.7cm}\left(256s^5-1600s^3t-27s^2t^4+2250st^2+108t^5+3125\right)^3.
\end{array}\label{kspecial}
\end{equation}
Finally, formulas (\ref{invariantk}), (\ref{numerator5}), (\ref{denominator5}), (\ref{i12}) yield
\begin{equation}
\begin{array}{l}
\displaystyle \ell(s,t):=L(f_{s,t})=-\frac{225\cdot 10!}{4}(125-3st^2)\left(19200 s^6t^2-\right.\\
\vspace{-0.3cm}\\
\hspace{1.6cm}\left.160000 s^4t^3-1120 s^3t^6+440000 s^2t^4+3600 st^7+\right.\\
\vspace{-0.3cm}\\
\hspace{1.6cm}\left.27 t^{10}-400000 t^5\right)/\left(256 s^5-1600s^3t-27s^2t^4+\right.\\
\vspace{-0.3cm}\\
\hspace{1.6cm}\left.2250st^2+108t^5+3125\right)^2.
\end{array}\label{lspecial}
\end{equation}

Now  Proposition \ref{equivarbitr} and Theorem \ref{solinv} imply the following result.

\begin{theorem}\label{specialequivsol} \sl Two germs ${\mathcal V}_{s_1,t_1}$ and ${\mathcal V}_{s_2,t_2}$ are biholomorphically equivalent if and only if
\begin{equation}
j(s_1,t_1)=j(s_2,t_2),\quad k(s_1,t_1)=k(s_2,t_2),\quad \ell(s_1,t_1)=\ell(s_2,t_2),\label{system}
\end{equation}
where the functions $j$, $k$, $\ell$ are given by {\rm (\ref{jspecial})}, {\rm (\ref{kspecial})}, {\rm (\ref{lspecial})}.
\end{theorem}
In the next section we will give an alternative proof of Theorem \ref{specialequivsol} based on computing the values of certain invariants for some binary sextic ${\bf f}_{s,t}$ arising from $f_{s,t}$.

\section{Associated binary forms}\label{assocforms}
\setcounter{equation}{0}

Let $Q(z,w)$ be a form in $X^n$ with $n\ge 3$ (see (\ref{setx})), and ${\mathcal V}$ the germ of the plane curve $\{Q=0\}$ at the origin. Then the germ of $Q$ lies in the {\it Jacobian ideal}\, ${\mathcal J}(Q)$ in ${\mathcal O}_2$, which is the ideal generated by the germs of all first-order partial derivatives of $Q$. Hence the moduli algebra ${\mathcal A}({\mathcal V})$ of ${\mathcal V}$ coincides with the {\it Milnor algebra}\, ${\mathcal O}_2/{\mathcal J}(Q)$. The germs of the partial derivatives of $Q$ form a regular sequence in ${\mathcal O}_2$, and therefore ${\mathcal A}({\mathcal V})$ is a complete intersection ring, which implies that ${\mathcal A}({\mathcal V})$ is Gorenstein (see \cite{B}). Recall that a local commutative associative algebra of finite vector space dimension greater than 1 is a {\it Gorenstein ring}\, if and only if for the annihilator $\Ann({\mathfrak m}):=\{u\in{\mathfrak m}: u\cdot{\mathfrak m}=0\}$ of its maximal ideal ${\mathfrak m}$ we have $\dim_{\CC}\Ann({\mathfrak m})=1$ (see, e.g. \cite{H}). In addition, ${\mathcal A}({\mathcal V})$ is {\it (non-negatively) graded}, that is, one has ${\mathcal A}({\mathcal V})=\oplus_{j\ge0}{\mathcal A}_{j}$, ${\mathcal A}_{j}{\mathcal A}_{k}\subset {\mathcal A}_{j+k}$, where ${\mathcal A}_{j}$ are linear subspaces of ${\mathcal A}({\mathcal V})$, with ${\mathcal A}_0\simeq\CC$.

Let ${\mathfrak m}({\mathcal V})$ be the maximal ideal of ${\mathcal A}({\mathcal V})$. Define the {\it exponential map}\, $\exp: {\mathfrak m}({\mathcal V})\ra {\bf 1}+{\mathfrak m}({\mathcal V})$ as follows:
$$
\displaystyle\exp(u):={\bf 1}+\sum_{k=1}^{\infty}\frac{1}{k!}u^k,\quad u\in{\mathfrak m}({\mathcal V}),
$$
where ${\bf 1}$ is the unit of ${\mathcal A}({\mathcal V})$. By Nakayama's lemma, ${\mathfrak m}({\mathcal V})$ is a nilpotent algebra, and therefore the above sum is in fact finite, with the highest-order term corresponding to $k=\nu$, where $\nu\ge 2$ is the {\it nil-index}\, of ${\mathfrak m}({\mathcal V})$ (i.e. the largest of all integers $\mu$ for which ${\mathfrak m}({\mathcal V})^{\mu}\ne 0$). Thus, the exponential map is a polynomial transformation.

Let $\pi$ be a projection on ${\mathfrak m}({\mathcal V})$ with range $\Ann({\mathfrak m}({\mathcal V}))$ (we call such projections {\it admissible}). Define
$$
S_{\pi}:=\left\{u\in{\mathfrak m}({\mathcal V}):\pi(\exp(u))=0\right\},
$$
where $\pi$ is extended to all of ${\mathcal A}({\mathcal V})$ by the condition $\pi({\bf 1})=0$. Observe that $S_{\pi}$ is an algebraic hypersurface in ${\mathfrak m}({\mathcal V})$ passing through the origin. Indeed, let ${\mathcal K}:=\ker\pi$ and $m:=\dim_{\CC}{\mathcal K}=\dim_{\CC}{\mathfrak m}({\mathcal V})-1$ (note that $m\ge 2$). Choose coordinates $\zeta_0,\zeta=(\zeta_1,\dots,\zeta_m)$ in ${\mathfrak m}({\mathcal V})$ such that ${\mathcal K}=\{\zeta_0=0\}$ and $\Ann({\mathfrak m}({\mathcal V}))=\{\zeta=0\}$. In these coordinates $S_{\pi}$ is given by a polynomial equation $\zeta_0=P_{\pi}(\zeta)$, where $P_{\pi}$ has neither constant nor linear term and $\deg P_{\pi}=\nu$. 

In Theorem 2.11 in \cite{FIKK} we obtained a criterion for two complex graded Gorenstein algebras of finite vector space dimension greater than 1 to be isomorphic (see also \cite{FK}, \cite{I} for results on algebras over arbitrary fields of characteristic zero). Applying this criterion to the moduli algebras of homogeneous plane curve germs we see that two moduli algebras ${\mathcal A}({\mathcal V_1})$, ${\mathcal A}({\mathcal V_2})$ are isomorphic if and only if for some (hence for any) admissible projections $\pi_1$, $\pi_2$ on ${\mathfrak m}({\mathcal V_1})$, ${\mathfrak m}({\mathcal V_2})$, respectively, and for some (hence for any) choice of coordinates in these algebras as described above, the polynomials $P_{\pi_1}$, $P_{\pi_2}$ are {\it linearly equivalent up to scale}, that is, there exist $c\in\CC^{*}$ and $C\in\GL(m,\CC)$ with 
\begin{equation}
c{\cdot}P_{\pi_1}(\zeta)=P_{\pi_2}(C\zeta).\label{mainids}
\end{equation}
Thus, the germs ${\mathcal V_1}$ and ${\mathcal V_2}$ are biholomorphically equivalent if and only if some (hence any) polynomials $P_{\pi_1}$, $P_{\pi_2}$ arising from ${\mathcal A}({\mathcal V_1})$, ${\mathcal A}({\mathcal V_2})$ are linearly equivalent up to scale (see Theorem 3.3 in \cite{FIKK}). We note that the linear equivalence of the polynomials $P_{\pi_1}$, $P_{\pi_2}$ up to scale takes place if and only if the hypersurfaces $S_{\pi_1}$, $S_{\pi_2}$ are {\it affinely}\, equivalent.

Clearly, identity (\ref{mainids}) can be rewritten as
$$
c{\cdot}P^{[l]}_{\pi_1}(\zeta)=P^{[l]}_{\pi_2}(C\zeta),\quad l=2,\dots\nu,
$$
where $P_{\pi}^{[l]}$ denotes the homogeneous component of order $l$ of $P_{\pi}$. We will now concentrate on the highest-order term $P_{\pi}^{[\nu]}$ of $P_{\pi}$. If $Q$ has degree $n$, it follows that $\nu=2(n-2)$ since $\Ann({\mathfrak m}({\mathcal V}))$ is spanned by the element of ${\mathfrak m}({\mathcal V})$ represented by the germ of the Hessian of $Q$ (see \cite{Sa}). We now choose the coordinates $\zeta_1,\dots,\zeta_m$ in such a way that for the corresponding basis $e_1,\dots,e_m$ in ${\mathcal K}$ the vectors $e_3,\dots,e_m$ span ${\mathfrak m}({\mathcal V})^2$ and the vectors $e_1$, $e_2$ span a complement to ${\mathfrak m}({\mathcal V})^2$ in ${\mathfrak m}({\mathcal V})$ (observe that $\dim_{\CC}{\mathfrak m}({\mathcal V})/{\mathfrak m}({\mathcal V})^2=2$). It is then clear that $P_{\pi}^{[\nu]}$ depends only on the variables $\zeta_1$, $\zeta_2$, that is, $P_{\pi}^{[\nu]}$ is a binary form of degree $2(n-2)$ on the subspace spanned by $e_1$, $e_2$. We note in passing that in fact $P_{\pi}$ can be introduced in a coordinate-free way as a polynomial on ${\mathfrak m}({\mathcal V})/\Ann({\mathfrak m}({\mathcal V}))$, in which case $P_{\pi}^{[\nu]}$ gives rise to a form on ${\mathfrak m}({\mathcal V})/{\mathfrak m}({\mathcal V})^2$.

For convenience, we will now make a canonical choice of $\pi$ and the vectors $e_1$, $e_2$. Let ${\mathcal M}$ (resp., ${\mathcal M}'$) be the collection of germs at the origin of all monomials $p$ in the variables $z$, $w$ with $1\le\deg p< 2(n-2)$ (resp., with $2\le\deg p\le 2(n-2)$), and define $\pi$ to be the admissible projection on ${\mathfrak m}({\mathcal V})$ whose kernel coincides with the span of the image of ${\mathcal M}$ in ${\mathcal A}({\mathcal V})$ under the factorization map ${\mathcal P}:{\mathcal O}_2\ra {\mathcal O}_2/{\mathcal J}(Q)={\mathcal A}({\mathcal V})$. Further, let $e_1$, $e_2$ be the images of the germs of the coordinate functions $z$, $w$, respectively, under ${\mathcal P}$, and observe that ${\mathcal P}\left({\mathcal M}'\right)$ spans ${\mathfrak m}({\mathcal V})^2$. With this choice of $\pi$, $e_1$, $e_2$, the form $P_{\pi}^{[\nu]}$ is defined up to a non-zero scalar factor, the only freedom being the choice of the coordinate $\zeta_0$ in $\Ann({\mathfrak m}({\mathcal V}))$. We denote the collection of mutually proportional binary forms of degree $2(n-2)$ obtained in this way by $\hat Q$ and say that any element ${\bf Q}\in\hat Q$ is a {\it binary form associated to $Q$}. Note that the above construction works in any number of variables, but for the purposes of this paper we restrict our considerations to the case of binary forms. We also remark that $\hat Q$ was first introduced in \cite{Ea} in slightly different terms. 

For an absolute classical invariant ${\mathcal I}$ of binary forms of degree $2(n-2)$, the value ${\mathcal I}({\bf Q})$ does not depend on the choice of the form ${\bf Q}$ associated to $Q$. Furthermore, the results of \cite{FIKK} discussed above imply that if $Q_1$, $Q_2$ are linearly equivalent binary forms, then one has ${\mathcal I}({\bf Q}_1)={\mathcal I}({\bf Q}_2)$ for any ${\bf Q}_1$, ${\bf Q}_2$ associated to $Q_1$, $Q_2$, respectively (cf. Theorem 2.1 in \cite{Ea}). Observe also that ${\mathcal I}({\bf Q})$ is rational when regarded as a function of $Q$, with ${\bf Q}\in\hat Q$ and $Q\in X^n$.

Let ${\mathcal R}^n$ denote the collection of all invariant rational functions on $X^n$ obtained in this way. Further, let $\check{\mathcal I}^n$ be the algebra of the restrictions to $X^n$ of {\it all}\, absolute invariants of forms of degree $n$ on $\CC^2$. Note that ${\mathcal R}^n$ lies in $\check{\mathcal I}^n$ (see Proposition 1 in \cite{DC}). We propose the following conjecture.

\begin{conjecture}\label{conj} \rm ${\mathcal R}^n=\check{\mathcal I}^n$.
\end{conjecture}

\noindent Since every element of $\check{\mathcal I}^n$ can be represented as a ratio of two elements of ${\mathcal I}^n$ (see Proposition 6.2 in \cite{M}), Conjecture \ref{conj} is equivalent to the statement ${\mathcal I}^n\subset {\mathcal R}^n$.

If Conjecture \ref{conj} is confirmed, it would provide a procedure for extracting a set of invariants of homogeneous plane curve singularities from their moduli algebras that solves the biholomorphic equivalence problem for such singularities (see Proposition \ref{equivarbitr}). This would be a step towards understanding how a germ ${\mathcal V}$ can be explicitly recovered from its moduli algebra ${\mathcal A}({\mathcal V})$ in general.

As we will see below, the families ${\mathcal V}_t$ and ${\mathcal V}_{s,t}$ considered in Sections \ref{proofs1} and \ref{n5}, respectively, supply evidence in favor of Conjecture \ref{conj} (see Remarks \ref{rem1} and \ref{rem2}). Specifically, we will give alternative proofs of Theorems \ref{main1} and \ref{specialequivsol} based on applying classical invariant theory to relevant associated forms rather than the forms $f_t$, $f_{s,t}$ as was done above. Before doing this, however, we will discuss the easier case of binary quartics.

Every non-zero square-free binary quartic is linearly equivalent to a quartic of the form
$$
q_t(z,w):=z^4+tz^2w^2+w^4,\quad t\ne\pm 2
$$ 
(see pp. 277--279 in \cite{El}). Any binary form associated to $q_t$ is again a quartic and is proportional to
$$
{\bf q}_t(\zeta_1,\zeta_2):=t\zeta_1^4-12\zeta_1^2\zeta_2^2+t\zeta_2^4
$$
(see \cite{Ea}). For $t\ne 0,\pm 6$ the quartic ${\bf q}_t$ is square-free, in which case the original quartic $q_t$ is associated to ${\bf q}_t$, and it is reasonable to say that for $t\ne 0, \pm 2, \pm 6$ the quartics $q_t$ and ${\bf q}_t$ are {\it dual}\, to each other. 

The algebra of classical invariants of binary quartics is generated by certain invariants ${\mathsf I}_2$ and ${\mathsf I}_3$, where, as usual, the subscripts indicate the degrees (see, e.g. pp. 101--102 in \cite{El}). For a quartic of the form
$$
Q(z,w)=a_4z^4+6a_2z^2w^2+a_0w^4
$$
the values of the invariants ${\mathsf I}_2$ and ${\mathsf I}_3$ are computed as follows: 
$$
\begin{array}{l}
\displaystyle {\mathsf I}_2(Q)=a_0a_4+3a_2^2=\frac{(Q,Q)^{(4)}}{1152},\quad {\mathsf I}_3(Q)=a_0a_2a_4-a_2^3,
\end{array}
$$
(see (\ref{transvectant})), and $\Delta(Q)={\mathsf I}_2(Q)^3-27\,{\mathsf I}_3(Q)^2$. Define an absolute invariant of binary quartics as
$$
{\mathsf J}:=\frac{{\mathsf I}_2^3}{\Delta}=\frac{J}{1152^3}
$$
(see (\ref{invariantj})). The restriction ${\mathsf J}|_{X^4}$ generates the algebra ${\mathcal I}^4$, and we have
$$
{\mathsf J}(q_t)=\frac{(t^2+12)^3}{108(t^2-4)^2}.
$$

Consider another absolute invariant of binary quartics
\begin{equation}
{\mathsf K}:=\frac{{\mathsf I}_2^3}{27{\mathsf I}_3^2}.\label{invarrmk}
\end{equation}
Then one obtains ${\mathsf K}({\bf q}_t)={\mathsf J}(q_t)$, and therefore ${\mathsf K}({\bf Q})={\mathsf J}(Q)$ for any $Q\in X^4$ and any ${\bf Q}\in\hat Q$. Thus, the absolute invariant ${\mathsf K}$ evaluated for associated quartics yields a generator of ${\mathcal I}^4$, which agrees with Conjecture \ref{conj}.

We will now give alternative proofs of Theorems \ref{main1} and \ref{specialequivsol}.
\vspace{0.3cm}

\noindent{\bf Proof of Theorem \ref{main1}:} We need to obtain just the \lq\lq only if\rq\rq\, implication. Assume first that $n\ge 5$ and find the value of the absolute invariant $M$ defined in (\ref{invariantm}) for forms associated to $f_t$. Any such form has degree $2(n-2)$ and is proportional to
$$
\begin{array}{l}
\displaystyle{\mathbf f}_t(\zeta_1,\zeta_2):=\sum_{j=n-1}^{2(n-2)}\left(\begin{array}{c}2(n-2)\\ j\end{array}\right)\left(\frac{(1-n)t}{n}\right)^{j+2-n}\zeta_1^j\zeta_2^{2(n-2)-j}+\\
\vspace{-0.1cm}\\
\displaystyle\hspace{2.3cm}\frac{(n-1)t^2}{n^2}\sum_{j=n-1}^{2(n-2)}\left(\begin{array}{c}2(n-2)\\j\end{array}\right)\left(\frac{(1-n)t}{n}\right)^{2(n-2)-j}\zeta_1^{2(n-2)-j}\zeta_2^j+\\
\vspace{-0.1cm}\\
\displaystyle\hspace{2.3cm}\left(\begin{array}{c}2(n-2)\\n-2\end{array}\right)\zeta_1^{n-2}\zeta_2^{n-2}.
\end{array}
$$
A straightforward computation yields
\begin{equation}
\begin{array}{l}
({\mathbf f}_t,{\mathbf f}_t)^{(2(n-2))}=((2(n-2))!)^2 \left(\begin{array}{c}2(n-2)\\n-2\end{array}\right)\Delta(f_t),\\
\vspace{-0.1cm}\\
({\mathbf H}({\mathbf f}_t),{\mathbf H}({\mathbf f}_t))^{(2(2n-6))}=\Delta(f_t)^2(\rho\Delta(f_t)+\sigma)
\end{array}\label{invark}
\end{equation}
for some $\rho,\sigma\in\CC$ with $\rho\ne 0$. The expressions in (\ref{invark}) imply
\begin{equation}
M({\mathbf f}_t)=\rho'\Delta(f_t)+\sigma'=\rho''t^n+\sigma''\label{kfinal}
\end{equation}
for some $\rho',\sigma',\rho'',\sigma''\in\CC$ with $\rho',\rho''\ne 0$, which leads to the desired result for $n\ge 5$.

Next, as shown above, for $n=4$ we have $1152^3\,{\mathsf K}({\mathbf f}_t)=J(f_t)$, where ${\mathsf K}$ is the absolute invariant of binary quartics defined in (\ref{invarrmk}). Formula (\ref{jphit}) now completes the proof.\qed

\begin{remark}\label{rem1} \rm It is clear from (\ref{kfinal}) that for $n\ge 5$ and suitable $a,b\in\CC$ the absolute invariant
$$
M'(Q):=\frac{\Bigl[(Q,Q)^{(2(n-2))}\Bigr]^2}{a({\mathbf H}(Q),{\mathbf H}(Q))^{(2(4n-6))}+b\Bigl[(Q,Q)^{(2(n-2))}\Bigr]^2}
$$
of forms of degree $2(n-2)$ has the property $M'({\mathbf f}_t)=J(f_t)$ for all $t$. Thus, one can recover the value of the absolute invariant $J$ for the form $f_t$ by evaluating a certain absolute invariant of binary forms of degree $2(n-2)$ for the associated form ${\bf f}_t$, which agrees with Conjecture \ref{conj}.
\end{remark}

\noindent{\bf Proof of Theorem \ref{specialequivsol}:} Any form associated to the quintic $f_{s,t}$ is proportional to the following binary sextic:
$$
\hspace{-0.2cm}\begin{array}{l}
{\bf f}_{s,t}(\zeta_1,\zeta_2):=(160s^3-300 st-27t^4)\zeta_1^6+(-1200s^2+81st^3+1125t)\zeta_1^5\zeta_2+\\
\vspace{-0.3cm}\\
\hspace{2.4cm}(-270s^2t^2+3750s+675t^3)\zeta_1^4\zeta_2^2+(480s^3t-1650st^2-\\
\vspace{-0.3cm}\\
\hspace{2.4cm}6250)\zeta_1^3\zeta_2^3+(-480s^4+2100s^2t-1125t^2)\zeta_1^2\zeta_2^4+(240s^3+\\
\vspace{-0.3cm}\\
\hspace{2.4cm}27s^2t^3-825st-108t^4)\zeta_1\zeta_2^5+(-6s^3t^2-50s^2+24st^3+125t)\zeta_2^6.
\end{array}
$$
Further, the algebra of classical invariants of binary sextics is generated by certain invariants of degrees 2, 4, 6, 10, 15 (see, e.g. \cite{Sy}, pp. 322--325 in \cite{El}, \cite{Ea}). We use the invariants ${\mathbf I}_2$, ${\mathbf I}_4$, ${\mathbf I}_6$, ${\mathbf I}_{10}$, ${\mathbf I}_{15}$ (with ${\mathbf I}_2(Q):=(Q,Q)^{(6)}/(6!)^2$) utilized in \cite{Ea}, where, as before, the subscripts indicate the degrees. 

Consider two absolute invariants of binary sextics
$$
{\mathbf J}:=\frac{3}{5}\frac{{\mathbf I}_2^2}{{\mathbf I}_2^2-2{\mathbf I}_4},\quad {\mathbf K}:=759375\frac{{\mathbf I}_{10}^2}{({\mathbf I}_2^2-2{\mathbf I}_4)^5}.
$$
In \cite{Ea} the values ${\mathbf J}({\bf f}_{s,t})$ and ${\mathbf K}({\bf f}_{s,t})$ were computed as follows:
\begin{equation}
\begin{array}{l}
\displaystyle {\mathbf j}(s,t):={\mathbf J}({\bf f}_{s,t})=(125-3st^2)^2/\\
\vspace{-0.3cm}\\
\hspace{1.65cm}(256s^5-1600s^3t-27s^2t^4+2250st^2+108t^5+3125),\\
\vspace{-0.1cm}\\
\displaystyle {\mathbf k}(s,t):={\mathbf K}({\bf f}_{s,t})=F(s,t)^2/\\
\vspace{-0.3cm}\\
\hspace{1.75cm}\left(256s^5-1600s^3t-27s^2t^4+2250st^2+108t^5+3125\right)^3,
\end{array}\label{bfjk}
\end{equation}
where
$$
\begin{array}{l}
F(s,t):=163200s^6 t^2+14800000s^5-2100000s^4t^3+5400s^3t^6-\\
\vspace{-0.3cm}\\
\hspace{1.9cm}92500000s^3 t+7425000s^2 t^4-52650st^7+116250000st^2+\\
\vspace{-0.3cm}\\
\hspace{1.9cm}729 t^{10}-4556250 t^5+312500000.
\end{array}
$$
Now, from (\ref{numerator5}), (\ref{denominator5}), (\ref{i12}) we obtain
\begin{equation}
\begin{array}{l}
\displaystyle F(s,t)=-27\cdot 10^{10}\, I_{12}(f_{s,t})+\frac{115625}{4608\cdot 10!}\Delta(f_{s,t})(f_{s,t}^2,f_{s,t}^2)^{(10)}+\\
\vspace{-0.3cm}\\
\displaystyle \hspace{1.8cm}\frac{5}{2(19200\cdot 10!)^3}\Bigl((f_{s,t}^2,f_{s,t}^2)^{(10)}\Bigr)^3,
\end{array}\label{expressF}
\end{equation}
and formulas (\ref{invariantj}), (\ref{invariantk}), (\ref{denominator5}), (\ref{jspecial}), (\ref{bfjk}), (\ref{expressF}) yield
\begin{equation}
\begin{array}{l}
\displaystyle{\mathbf j}=\frac{1}{5(1440000\cdot 10!)^2}  j,\\
\vspace{-0.1cm}\\
{\mathbf k}= 2^{20}3^65^5\, k+\hbox{const}\, \ell  + \hbox{const}\, j \ell+\\
\vspace{-0.3cm}\\
\hspace{0.8cm}\hbox{const}\, j^3+\hbox{const}\, j^2+\hbox{const}\, j.
\end{array}\label{connection1}
\end{equation} 

Consider a third absolute invariant of binary sextics 
$$
{\mathbf L}:=675\frac{{\mathbf I}_2{\mathbf I}_{10}}{({\mathbf I}_2^2-2{\mathbf I}_4)^3}.
$$ 
We have\footnote{We thank M. Eastwood for his help with computing the value ${\mathbf L}({\bf f}_{s,t})$. The calculations have been performed using Maple.}
\begin{equation}
\begin{array}{l}
\displaystyle {\mathbf l}(s,t):={\mathbf L}({\bf f}_{s,t})=(125-3st^2)F(s,t)/\\
\vspace{-0.3cm}\\
\hspace{1.65cm}(256s^5-1600s^3t-27s^2t^4+2250st^2+108t^5+3125)^2.
\end{array}\label{bfl}
\end{equation}
Formulas (\ref{invariantj}), (\ref{invariantk}), (\ref{numerator5}), (\ref{denominator5}), (\ref{expressF}), (\ref{bfl}) yield
\begin{equation}
{\mathbf l}= -\frac{12}{25\cdot 10!}\,\ell+\hbox{const}\, j^2+\hbox{const}\, j.\label{connection2}
\end{equation} 
It follows from identities (\ref{connection1}), (\ref{connection2}) that system (\ref{system}) is equivalent to the system
$$
{\mathbf j}(s_1,t_1)={\mathbf j}(s_2,t_2),\quad {\mathbf k}(s_1,t_1)={\mathbf k}(s_2,t_2),\quad {\mathbf l}(s_1,t_1)={\mathbf l}(s_2,t_2),
$$
which completes the proof.\qed

\begin{remark}\label{rem2}\rm In the above proof we showed that one can recover the values of the absolute invariants $J$, $K$, $L$ for the quintic $f_{s,t}$ by evaluating certain absolute invariants of binary sextics for the associated form ${\bf f}_{s,t}$, which agrees with Conjecture \ref{conj}.
\end{remark}

\end{document}